\documentclass{ws-procs9x6}

\newcommand\slint{{\int\!\!\!\!\!\!{-}}}
\newcommand\tslint{{\int\!\!\!\!\!{-}}}
\newcommand\dslash{d\!\!{}^{\text{\rm--}}\!}
\newcommand\tr{\operatorname{tr}}
\newcommand\Tr{\operatorname{Tr}}
\newcommand\TR{\operatorname{TR}}
\newcommand\res{\operatorname{res}}

\begin{document}

\title{Remarks on nonlocal \\ trace expansion coefficients 
}

\author{Gerd Grubb}

\address{Mathematics Department, 
University of Copenhagen,\\
Universitetsparken 5, \\ 
2100 Copenhagen, Denmark\\ 
E-mail: {\tt grubb@math.ku.dk}}

\maketitle{}

\abstracts{In a recent work, Paycha and Scott establish formulas for
all the Laurent coefficients of $\Tr(AP^{-s})$ at the possible poles. In
particular, they show a formula for the zero'th coefficient at $s=0$,
in terms of two functions generalizing, respectively, the
Kontsevich-Vishik canonical trace density, and the Wodzicki-Guillemin
noncommutative residue density of an associated operator. The purpose
of this note is to 
provide a proof of that formula relying entirely on resolvent
techniques (for the sake of possible generalizations to situations
where powers are not an easy tool). --- We also give some corrections to
transition formulas used in our earlier works.  }


\section{Introduction}

In an interesting new work [PS], Sylvie Paycha and Simon Scott have
obtained formulas for all the coefficients in Laurent expansions of
zeta functions $\zeta (A,P,s)=\Tr(AP^{-s})$ around the poles, in
terms of combinations 
of  finite part integrals and residue type integrals, of associated
logarithmic symbols. We consider classical pseudodifferential
operators ($\psi $do's) $A$ and $P$ of order $\sigma \in\mathbb R$
resp.\ $m\in\mathbb R_+$ acting in a Hermitian vector bundle $E$ over a 
closed $n$-dimensional manifold $X$, $P$ being elliptic with
principal symbol eigenvalues in $\mathbb C\setminus \mathbb R_-$.
 The basic formula is the following formula for
$C_0(A,P)$, where $C_0(A,P)-\Tr(A\Pi _0(P))$ is the
regular value of $\zeta (A,P,s)$ at $s=0$:
\begin{equation} 
C_0(A,P)=\int_X\bigl(\TR_x(A)-\tfrac1m\operatorname{res}_{x,0}(A\log
P)\bigr)\,dx.
\label{0.1}
\end{equation}
The integrand is defined in a local coordinate system by:
\begin{equation} 
\TR_x(A)=\slint
\tr a(x,\xi )\,\dslash\xi,\;\operatorname{res}_{x,0}(A\log P)
=\int_{|\xi |=1}\tr r_{-n,0}(x,\xi )\dslash S(\xi ),
\label{0.2}
\end{equation}
where $\tslint  a(x,\xi )\,\dslash\xi$ is a finite part integral, and
$r$ is the symbol of $R=A\log P$, of the form $r(x,\xi
)\sim \sum_{j\ge 0,\, l=0,1}r_{\sigma -j,l}(x,\xi )(\log[\xi ])^l$;
$r_{\sigma -j,l}$ homogeneous in $\xi $ of degree $\sigma -j$
for $|\xi |\ge 1$, $[\xi ]=|\xi |$ for $|\xi |\ge 1$.
(Here $r_{-n,0}$ is set equal to 0 when $\sigma -j$ does not hit
$-n$.) Moreover, the expression $\bigl(\TR_x(A)-\operatorname{res}_{x,0}(A\log
P)\bigr)\,dx$ has an invariant meaning as a density on $X$, although
its two terms individually do not so in general. (In these formulas
we use the conventions $\dslash \xi =(2\pi )^{-n}d\xi $,
$\dslash S(\xi) =(2\pi )^{-n}dS(\xi )$, where $dS(\xi )$ indicates the usual
surface measure on the unit sphere.)

Formula (\ref{0.1}) generalizes the formula $C_0(A,P)=\TR A$ (the
canonical trace), which holds
in particular cases, cf.\ Kontesevich and Vishik [KV], Lesch [L], Grubb
[G4]. The general formula is shown in [PS] by use of
holomorphic families of 
$\psi $do's (depending holomorphically on their complex order $z$);
in particular complex powers of $P$.
The purpose of this note is to derive it by methods relying on the
knowledge of the resolvent $(P-\lambda )^{-1}$.
This is meant to facilitate
generalizations to manifolds with boundary, where powers of operators
are not an easy tool.

--- We take the opportunity here to correct, in the appendix, some
inaccuracies in earlier 
papers, mainly concerning the relations between expansion coefficients in
resolvent traces and zeta functions. 

\section{Preliminaries}

Recall the expansion formulas for the resolvent kernel and trace, worked out in
local coordinates in Grubb and Seeley [GS1], when $N>\frac
{\sigma +n}m$:
\begin{equation}\aligned
K\bigl( A(P-\lambda )^{-N},& x,x\bigr)
\sim
\sum_{ j\in \mathbb N  } \tilde c^{(N)}_{ j}(x)(-\lambda ) ^{\frac{\sigma
+n -j}m-N}\\
&+  
\sum_{k\in \mathbb N}\bigl( \tilde c^{(N)\prime}_{ k}(x)\log (-\lambda )
+\tilde c^{(N)\prime\prime}_{ k}(x)\bigr)(-\lambda ) ^{{ -k}-N},\\
\Tr\bigl( A(P-\lambda )^{-N}&\bigr)
\sim
\sum_{ j\in \mathbb N  } \tilde c^{(N)}_{ j}(-\lambda ) ^{\frac{\sigma
+n -j}m-N}\\
&+  
\sum_{k\in \mathbb N}\bigl( \tilde c^{(N)\prime}_{ k}\log (-\lambda )
+\tilde c^{(N)\prime\prime}_{ k}\bigr)(-\lambda ) ^{{ -k}-N},
\endaligned\label{1.1}\end{equation}
for $\lambda \to\infty $ on rays in 
a sector around $\mathbb R_- $. The second formula is deduced from the
first one by integrating the fiber trace in $x$. We denote
$\{0,1,2,\dots\}=\mathbb N$.
(More precisely, [GS1] covers the cases where $m$ is integer; the
noninteger cases are included in Loya [Lo], Grubb and Hansen [GH].)

The $\tilde c^{(N)\prime}_k(x)$ and $\tilde c^{(N)\prime}_k$ 
vanish when $\sigma +n+mk\notin \mathbb N$;
this holds for all $k$ when $m$ is integer and $\sigma $ is noninteger.
When $\sigma +n+mk=j\in \mathbb N$, $\tilde c^{(N)}_{ j}(x)$ and
$\tilde c^{(N)\prime\prime}_{ k}(x)$ are both coefficients of the
power $(-\lambda )^{-k-N}$; their individual values depend on the
localization used (as worked out in detail in [G4]), and it is only
the sum $\tilde c^{(N)}_{ j}+\tilde c^{(N)\prime\prime}_{ k}$ that
has an invariant meaning.

The coefficients depend on $N$; when $N=1$, we drop the upper index $(N)$. 
We are particularly interested in the coefficient of $(-\lambda
)^{-N}$, for which we shall use the notation
\begin{equation}
\widetilde C^{(N)}_0(A,P)= \tilde c^{(N)}_{\sigma +n}+\tilde
c^{(N)\prime\prime}_0;
\label{1.2}
\end{equation}
here we have for convenience set
\begin{equation}\tilde c^{(N)}_{\sigma +n}=0\text{ when
}\sigma +n\notin \mathbb 
N.\label{1.3}\end{equation}

Recall that there is, equivalently to the last expansion in (\ref{1.1}), an
expansion formula for complex powers:
\begin{equation}
\Gamma (s)\Tr (AP^{-s})\sim
\sum_{ j\in \mathbb N } \frac{  c_{j}}{s+\frac{j-\sigma -n}m}-\frac{\Tr (A\Pi
_0(P))}s +  \sum_{k\in \mathbb N}\Bigl(\frac{
  c'_{k}}{(s+k)^2}+\frac{  c''_{k}}{s+ k}\Bigr) .\label{1.4}
\end{equation} 
This means that $\Gamma (s)\Tr (AP^{-s})$,
defined as a holomorphic function for $\operatorname{Re}s>\frac
{\sigma +n}m$,  
extends meromorphically to $\mathbb C$ with the pole structure 
indicated in the right hand side. Here $
\Pi _0(P)=\tfrac i{2\pi }\int_{|\lambda |=\varepsilon }(P-\lambda
)^{-1}\,d\lambda $ is the projection onto
the generalized nullspace of $P$ (on which $P^{-s}$ is taken to be
zero). Again, the $ c'_k$ vanish when $\sigma +n+mk\notin \mathbb N$.
We denote 
\begin{equation}
C_0(A,P)=c_{\sigma +n}+c''_0,\label{1.5}
\end{equation}
{\it the basic coefficient} (setting $c_{\sigma +n}$ equal to 0 when $\sigma
+n\notin \mathbb N$). 

The transition between (\ref{1.1}) and (\ref{1.4}) is accounted for e.g.\ in Grubb and Seeley [GS2],
Prop.\ 2.9, (3.21). 
The coefficient sets in (\ref{1.1}) and (\ref{1.4}) are derivable from one
another. The
coefficients $\tilde c^{(N)}_j$ and $  c_j$, resp.\ $\tilde c^{(N)\prime}_k$ and $  c'_k$, 
are proportional by universal nonzero
constants. This holds also for $\tilde c^{(N)\prime\prime}_k$ and $
c''_k$, when the $c'_k$ vanish. In general,
there are linear formulas for the transitions between $\{\tilde
c^{(N)\prime}_k,\tilde c^{(N)\prime\prime}_k\}$ and $\{c'_k,  c''_k\}$. 
(For $N=1$, [GS2], Cor.\ 2.10, would imply that
$\tilde c''_k$ and $c''_k$ are proportional in general, but in fact, 
the formulas for the
$\tilde a_{j,l}$ given there 
are only correct for $l=m_j$, whereas for
$l<m_j$ there is an effect from derivatives of the gamma function
that was overlooked.)
One has in particular that $\tilde c^{(N)\prime}_0=  c'_0$ for all
$N$.

Division of (\ref{1.4}) by $\Gamma (s)$ gives the pole
structure of $\zeta (A,P,s)$:
\begin{equation}
\zeta (A,P,s)=\Tr (AP^{-s})\sim
\sum_{ j\in \mathbb N} \frac{
c'''_{j}}{s+\frac{j-\sigma -n}m},
 \label{1.6}
\end{equation}
 where $c'''_j$ is proportional to $c_j $ if  $\tfrac{j-\sigma
-n}m\notin\mathbb N$, and $c'''_j$ is proportional to $c'_k$ if $\tfrac{j-\sigma -n}m=k\in\mathbb N$. 
 
One can study the Laurent series expansions of $\zeta (A,P,s)$ at the
poles by use of (\ref{1.1}). We now restrict the
attention to the possible pole at $s=0$.

Write the Laurent expansion at 0 as follows:\begin{equation}
\zeta (A,P,s)=  C_{-1}(A,P)s^{-1}+(C_0(A,P)-\Tr(A\Pi _0(P)))s^0+\sum _{l\ge
1}C_{l}(A,P)s^l.\label{1.7}
\end{equation}
It is known from  Wodzicki [W], Guillemin [Gu], that
\begin{equation}
C_{-1}(A,P)= {  c'_0}=\tfrac1m\operatorname{res}(A),\text{
independently of }P;\label{1.8}
\end{equation}
it vanishes if $\sigma +n\notin\mathbb N$ or the symbols have certain
parity properties. From
 [KV],  [L], [G4] we have  
that $C_0(A,P)=\TR A$ when $\sigma +n\notin \mathbb N$, and in certain
parity cases (given in [KV] for $n$ odd, [G4] for $n$ even, more
details at the end of Section 4). Also $C_1(A,P)$
is of interest, since the zeta determinant of $P$ satisfies\begin{equation}
\log\det P=-C_1(I,P)=C_0(\log P,P)\label{1.9}
\end{equation}
(cf.\ Okikiolu [O], [G4]); here it is useful to know that
expansions 
like (\ref{1.1}) but with higher powers of $\log(-\lambda )$ hold if $A$ is
log-polyhomogeneous, cf.\ 
[L] and [G4].

We have for general $N\ge 1$:

\begin{lemma}\label{Lemma 1.1} When $N>\frac{\sigma +n}m$ (so that
$A(P-\lambda 
)^{-N}$ is trace-class), then
\begin{equation}\aligned
\tilde c^{(N)\prime}_0&=  c'_0=\tfrac1m\res A,\\
\widetilde C^{(N)}_0(A,P)&=C_0(A,P)-\alpha _N  c'_0,\text{ where}\\
\alpha _N&=1+\tfrac12+\dots+\tfrac1{N-1}.
\endaligned\label{1.10}
\end{equation}

\end{lemma}

\begin{proof} Denote $N-1=M$, then 
\begin{equation}(P-\lambda )^{-N}=(P-\lambda )^{-M-1}=\tfrac1{M!}\partial _\lambda
^M(P-\lambda 
)^{-1}.\label{1.11}\end{equation} 
The transition from (\ref{1.1}) to information on $\zeta (A,P,s)$ is
obtained by use of the formula \begin{equation}
AP^{-s}=\tfrac {M!}{(s-M)\dots(s-1)}\tfrac i{2\pi }\int_{\mathcal
C}\lambda ^{M-s}\tfrac1{M!}\partial _\lambda  ^MA(P-\lambda
)^{-1}(I-\Pi _0(P))\,d\lambda ,\label{1.12}
\end{equation}
where $\mathcal C$ is a curve in $\mathbb C\setminus {\overline{\mathbb R}_-}$ around the
nonzero spectrum of $P$. Here we can take traces on both sides and
apply [GS2], Prop.\ 2.9, to  
\begin{equation}f(\lambda )=\Tr\bigl(A\tfrac1{M!}\partial _\lambda ^M(P-\lambda
)^{-1}(I-\Pi _0(P))\bigr),\end{equation}  
defining \begin{equation}
\varrho (s)=\tfrac i{2\pi }\int_{\mathcal
C}\lambda ^{-s}f(\lambda )\,d\lambda,\label{1.13}
\end{equation}
for $\operatorname{Re}s$ large, and extending
meromorphically. Then 
\begin{equation}
\zeta
(A,P,s)=\tfrac{M!}{(s-M)\dots(s-1)}
\varrho (s-M).\end{equation}
(Note that $\zeta (A,P,s)=\varrho (s)$ if $N=1$.)
The  expansion coefficients of
$f(\lambda )$ in powers and log-powers are universally
proportional 
to the pole coefficients of 
$\psi (s)=\frac{\pi }{\sin (\pi s)}\varrho (s)$ at simple and double poles, for each
index, as accounted for in [GS2], Prop.\ 2.9. 

When we apply this to
$\zeta (A,P,s)$, we must take the factors
$g_M(s)=\tfrac{M!}{(s-M)\dots(s-1)}$ and $\frac1\pi \sin({\pi (s-M)})$ into
account. We have\begin{equation}
\zeta
(A,P,s)=g_M(s)\tfrac1\pi \sin(\pi (s-M))\psi (s-M)
.\label{1.14}\end{equation}
By [GS2], Prop.\ 2.9, a pair of terms $a(-\lambda
)^{-M-1}\log(-\lambda )+b(-\lambda )^{-M-1}$ in the expansion of
$f(\lambda )$ carries over to the pair of terms $\frac
a{(s+M)^2}+\frac b{s+M}$ in 
the pole structure of $\psi (s)$, whereby \begin{equation}
\psi (s-M)=\frac a{s^2}+\frac bs+O(1),\text{ for }s\to 0.
\end{equation}
Now it is easily checked that \begin{equation}\aligned
\tfrac1\pi \sin(\pi (s-M))&=(-1)^M(s+cs^3+O(s^5)),\\
g_M(s)&=(-1)^M(1+(1+\tfrac12+\dots+\tfrac1M)s+O(s^2)),
\endaligned\end{equation}
for $s\to 0$. Then, with $\alpha _N$ defined in (\ref{1.10}), 
\begin{equation}
\aligned
\zeta
(A,P,s)&=(s+cs^3+O(s^5))(1+\alpha _{M+1}s+O(s^2))\Bigl(\frac
a{s^2}+\frac bs+O(1)\Bigr)\\
&=\frac a{s}+(b+\alpha _{M+1} a)+O(s),\endaligned\label{1.15}
\end{equation}
for $s\to 0$. For $f(\lambda )$ in (\ref{1.13}) we have this situation with
$a=\tilde c^{(N)\prime}_0$ and $b$ \linebreak$=\widetilde
C^{(N)}_0(A,P)-\Tr(A\Pi _0(P))$, so (\ref{1.15}) holds with these values.
In view of (\ref{1.7}), (\ref{1.8}), this shows
 (\ref{1.10}).  
\end{proof}

If one writes $f(\lambda )$ in the lemma as a sum $f_1(\lambda
)+f_2(\lambda )$, where $f_1$ has the sum over $j$ in (\ref{1.1}),
resp.\ $f_2$ has the sum over $k$ in (\ref{1.1}), as asymptotic
expansions, the lemma 
can be applied to $f_1$ and $f_2$ separately, relating their
coefficients to those of the poles of the corresponding functions
of $s$. 

\begin{remark}\label{Remark 1.2} Formula (\ref{1.10}) gives a correction to our earlier papers
[G1--5] and Grubb and Schrohe  [GSc1--2], where it was
taken for granted that 
$C_0(A,P)$ would equal $\tilde c^{(N)}_{\sigma +n}+\tilde
c^{(N)\prime\prime}_0$ for any $N$. Fortunately, the correction has
no consequence for the results in those papers, which were either
concerned with the value of $C_0(A,P)$ when $  c'_0=0$,
or its value {\it modulo local terms} ($  c'_0$ is local), or values of
combined expressions where $  c'_0$-contributions cancel out. More on
corrections in the appendix.

In [GS2], Cor.\ 2.10 was not used in the argumentation, which
was based directly on Prop.\ 2.9 and the primary knowledge of zeta
expansions.

\end{remark}

We shall now analyze $C_0(A,P)$ further, showing (\ref{0.1}) by resolvent
considerations.
Our proof is based on an explicit calculation of
one simple special case, together with the use of the trace defect formula\begin{equation}
C_0(A,P)-C_0(A,P')=-\tfrac 1m \operatorname{res}(A(\log P-\log
P')).\label{1.16} 
\end{equation}
This formula is well-known from considerations of complex powers of
$P$ ([O],  [KV], Melrose and Nistor [MN]), but can also be
derived directly  
from resolvent considerations [G6].

\section{The trace defect formula for general orders}

In [G6], the arguments for (\ref{1.16}) are given in detail in cases
where $m>\sigma +n$, whereas more general cases are briefly explained by
reference to  Remark 3.12 there. For completeness, we give the explanation
in detail here. This is a minor technical point that may be skipped
in a first reading. Denote as in [G6]
\begin{equation}
S_\lambda =A((P-\lambda )^{-1}-(P'-\lambda )^{-1} ),
\label{2.1}\end{equation}
where $P$ and $P'$ are of order $m>0$; then (cf.\ (\ref{1.11}))
\begin{equation}
A((P-\lambda )^{-N}-(P'-\lambda )^{-N} )=\tfrac{\partial _\lambda
^{N-1}}{(N-1)!}S_\lambda \equiv S^{(N)}_\lambda .\label{2.2
}\end{equation}
$S_\lambda $ and $S^{(N)}_\lambda $ have symbols $s(x,\xi ,\lambda )$
resp.\ 
$s^{(N)}(x,\xi ,\lambda )=\tfrac{\partial _\lambda
^{N-1}}{(N-1)!}s(x,\xi ,\lambda )$ in local trivializations. 

The difference of the two
expansions (\ref{1.1}) with $P$ resp.\ $P'$ inserted satisfies
\begin{equation}
\aligned
\Tr&\bigl( S^{(N)}_\lambda \bigr)
\sim
\sum_{ j\in \mathbb N  } \tilde s^{(N)}_{ j}(-\lambda ) ^{\frac{\sigma
+n -j}m-N}\\
&+  \tilde s^{(N)\prime\prime}_{0}(-\lambda )^{-N}
+\sum_{k\ge 1}\bigl( \tilde
s^{(N)\prime}_{ k}\log (-\lambda ) +\tilde s^{(N)\prime\prime}_{ 
k}\bigr)(-\lambda ) ^{{ -k}-N};
\endaligned\label{2.3}\end{equation}
in view of Lemma \ref{Lemma 1.1}, the contributions from
$\operatorname{res}A$ cancel out and the coefficient of $(-\lambda
)^{-N}$ equals $C_0(A,P)-C_0(A,P')$.

The symbol $s(x,\xi ,\lambda )$ is analyzed in [G6], Prop.\ 2.1.
For $s^{(N)}$, we conclude that the homogeneous terms
have at least $N+1$ factors of the form $(p_{m}-\lambda )^{-1}$ or
$(p'_{m}-\lambda )^{-1}$, hence
the strictly homogeneous version of the symbol of order  $\sigma -Nm-j$
satisfies\begin{equation}
|s^{(N)h}_{\sigma -Nm-j}(x,\xi ,\lambda )|\le c (|\xi |^m+|\lambda
|)^{-N-1}|\xi |^{\sigma +m-j},\label{2.4}\end{equation}being
integrable at $\xi =0$ for $j<n+\sigma +m$, $\lambda \ne 0$. Then the
kernel and trace 
of $S^{(N)}_\lambda $ have expansions
\begin{equation}
\aligned
K(S^{(N)}_\lambda ,x,x)
&=
\sum_{ j<\sigma + m+n  }  \tilde s^{(N)}_{ j}(x)(-\lambda ) ^{\frac{n
+\sigma -j}m  -N}+
O(|\lambda |^{-N-1+\varepsilon }),\\
\Tr S^{(N)}_\lambda 
&=
\sum_{ j<\sigma +m+n }  \tilde s^{(N)}_{ j}(-\lambda ) ^{\frac{n +\sigma -j}m  -N}+
O(|\lambda |^{-N-1+\varepsilon }),
\endaligned
\label{ 2.5
}\end{equation}
where the terms for $j<\sigma +m+n$ are calculated from the
strictly homogeneous symbols (for $\lambda \in\mathbb R_-$):
\begin{equation}
\int_{\mathbb R^n}s^{(N)h}_{\sigma -Nm-j}(x,\xi ,\lambda )\,\dslash\xi
=(-\lambda )^{\frac{n+\sigma -j}m-N}\int_{\mathbb R^n}s^{(N)h}_{\sigma
-Nm-j}(x,\eta ,-1)\,\dslash\eta ,\label{2.6}
\end{equation}
so that 
\begin{equation}
  \tilde s^{(N)}_j(x)=\int_{\mathbb R^n}s^{(N)h}_{\sigma -Nm-j}(x,\xi ,-1)\,\dslash\xi ,\quad 
  \tilde s^{(N)}_j=\int \tr \tilde s^{(N)}_j(x)\,dx .
\end{equation} 

When $\sigma +n\notin\mathbb N$, there is no term
with $(-\lambda )^{-N}$ in the expansion of  $\Tr S^{(N)}_\lambda $, and
(\ref{1.16}) holds trivially. 

When $\sigma +n\in\mathbb N$, the coefficient of
$(-\lambda )^{-N}$ in 
$\Tr S^{(N)}_\lambda $ equals  
\begin{equation}
C_0(A,P)-C_0(A,P')=  \tilde s^{(N)}_{\sigma +n}=\int\int_{\mathbb R^n}\tr
s^{(N)h}_{-Nm-n}(x,\xi 
,-1)\,\dslash\xi dx .\label{2.7}\end{equation} 
This demonstrates that the term is local, and gives a means to
calculate it (as indicated in [G6], Rem.\ 3.12): Note that
$s^{(N)h}_{ -Nm-n}(x,\xi ,\lambda )=\frac {\partial _\lambda
^{N-1}}{(N-1)!}s^h_{-m-n}(x,\xi ,\lambda )$. Since $s^h_{-m-n}$
satisfies (\ref{2.4}) with $N=1$, $j=\sigma +n$, it is integrable over
$\mathbb R^n$ when $\lambda \ne 0$. Here \begin{equation}
(-\lambda )^{-1}\int_{\mathbb R^{n}} s^h_{-m-n}(x,\xi
,-1)\,\dslash\xi =\int _{\mathbb R^{n}} s^h_{-m-n}(x,\xi
,\lambda )\,\dslash\xi \end{equation}
for $\lambda \in\mathbb R_-$. Moreover, the integral from (\ref{2.6})
satisfies  
\begin{equation}\aligned
\int _{\mathbb R^{n}}  s^{(N)h}_{-Nm-n}(x ,\xi
 ,\lambda )\,\dslash\xi  &=\tfrac {\partial _\lambda
^{N-1}}{(N-1)!}\int _{\mathbb R^{n}}  s^{h}_{-m-n}(x ,\xi
 ,\lambda )\,\dslash\xi  \\
&=\tfrac {\partial _\lambda
^{N-1}}{(N-1)!}\bigl[(-\lambda )^{-1}\int_{\mathbb R^{n}}  s^h_{-m-n}(x ,\xi
 ,-1)\,\dslash\xi  \bigr]\\
&=(-\lambda )^{-N}\int_{\mathbb R^{n}}  s^h_{-m-n}(x ,\xi
 ,-1)\,\dslash\xi  ,\endaligned
\end{equation} 
which implies
\begin{equation}
  \tilde s^{(N)}_{\sigma +n}(x)=\int_{\mathbb R^n}s^{(N)h}_{ -Nm-n}(x,\xi ,-1)\,\dslash\xi =\int_{\mathbb R^{n}}  s^h_{-m-n}(x ,\xi
 ,-1)\,\dslash\xi .\label{2.8}\end{equation}
The latter in turned into the residue integrand for
$-\frac1m\operatorname{res}(A(\log P-\log P'))$ by Lemmas 1.2 and 1.3
of [G6], as already done in
Section 2 there.

We conclude:

\begin{theorem}\label{Theorem 2.1}  Let $P$ and $P'$ be classical elliptic $\psi $do's of
order $m\in\mathbb R_+$ and such that the principal symbol has no eigenvalues on
$\mathbb R_-$, let $A$ be a 
classical $\psi $do of order $\sigma $, and let
$S_\lambda =A((P-\lambda )^{-1} 
-(P'-\lambda )^{-1})$ and 
$F=A(\log P-\log P')$ with symbols $s$ resp.\ $f$. 

Consider the case $\sigma +n\in\mathbb N$.
Then
\begin{equation}
C_0(A,P)-C_0(A,P')=\int_X \tr \tilde s_{\sigma +n}(x)\,dx
=-\tfrac1m\operatorname{res}(A(\log P-\log P'))
\label{2.9}\end{equation} 
where, for each $x$, in local coordinates,
\begin{equation}
\tilde s_{\sigma +n}(x)=\int_{\mathbb R^n}s^h_{-m-n}(x,\xi ,-1)\,\dslash\xi 
=-\tfrac 1m \int_{|\xi
|=1}f_{-n}(x,\xi )\, \dslash S(\xi ).
\label{2.10}\end{equation}

When $\sigma +n\notin \mathbb N$, the identities hold trivially (with zero
values everywhere).
\end{theorem} 

It follows moreover:

\begin{corollary}\label{Corollary 2.2} If, in Theorem {\rm \ref{Theorem 2.1}}, $P'$ is replaced
by an operator of a different order $m'>0$, then one has:
\begin{equation}
C_0(A,P)-C_0(A,P')=
-\operatorname{res}(A(\tfrac1m\log P-\tfrac1{m'}\log P')).
\label{2.13}\end{equation} 

\end{corollary}

\begin{proof} Let $P_0$ be an elliptic, selfadjoint positive $\psi $do
of order
$m$, and define $P_0^{m'/m}$ by spectral calculus; it is an elliptic,
selfadjoint positive $\psi $do of order $m'$. Then by the definition
of the zeta function, $
C_0(A,P_0)=C_0(A,P_0^{m'/m})$. Applications of Theorem \ref{Theorem 2.1} with
$P,P_0$ and with 
$P',P_0^{m'/m}$ give:\begin{equation}
\aligned
C_0&(A,P)-C_0(A,P')=(C_0(A,P)-C_0(A,P_0))\\
&\qquad+
(C_0(A,P_0)-C_0(A,P_0^{m'/m}))+(C_0(A,P_0^{m'/m})-C_0(A,P'))\\
&=
-\operatorname{res}(A(\tfrac1m\log P-\tfrac 1m\log P_0))
-\operatorname{res}(A(\tfrac1{m'}\log P_0^{m'/m}-\tfrac 1{m'}\log P'))
\\
&=-\operatorname{res}(A(\tfrac1m\log P-\tfrac 1{m'}\log P')),
\endaligned
\end{equation}
as was to be shown.  
\end{proof}

\section{A formula for the zero'th coefficient}

Our strategy for calculating $C_0(A,P)$ is to use (\ref{2.13}) in combination
with an exact calculation for a very special choice
$P_0$ of $P$, namely \begin{equation}
P_0=((-\Delta )^{m/2}+1)I_M\text{ with symbol }p_0=(|\xi
|^m+1)I_M,\label{3.1}
\end{equation} 
in suitable local coordinates; here $I_M$ is the $M\times M$ identity
matrix (understood in the following), $M=\dim E$, and $m$ is even.

Let $A$ be given, of order 
$\sigma \in\mathbb R$, then we take $m>\sigma +n$. 
Let $\Phi _j:E|_{U_j}\to V_j\times \mathbb C^M$, $j=1,\dots,J$, be an
atlas of trivializations with base maps $\kappa _j$ from $U_j\subset
X$ to $V_j\subset \mathbb R^n$, let $\{\psi _j\}_{1\le j\le J}$ be an
associated partition of unity (with $\psi _j\in C_0^\infty (U_j)$),
and let $\varphi _j\in C_0^\infty (U_j)$ with $\varphi _j=1$ on
$\operatorname{supp}\psi _j$. Then 
\begin{equation}
A=\sum_{1\le j\le J}\psi _jA= \sum_{1\le j\le J}\psi _jA\varphi
_j+\sum_{1\le j\le J}\psi _jA(1-\varphi _j),
\end{equation}
where the last sum is a $\psi $do of order $-\infty $; for this the
formula (\ref{0.1}) is well-known, since $C_0(B,P)=\Tr B$ when $B$ is
of order $<-n$. So it remains to treat each of the terms $\psi _jA\varphi
_j$. Consider e.g.\ $\psi _1A\varphi _1$. We could have
assumed from the start that $X$ was already covered by a family of
open subsets $U_{j0}\subset \subset U_j$. Thus it is no restriction to
assume that $\psi _1$ and $\varphi _1$ are supported in
$U_{10}\subset\subset U_1$, where  $U_{10}, U_2,\dots, U_J$ cover $X$.

Replace $U_j$ by $U'_j=U_j\setminus \overline U_{10}$ for $j\ge 2$,
and write $U_1=U'_1$, then $\{U'_j\}_{1\le j\le J}$ also covers $X$.
Let $\{\psi '_j\}_{1\le j\le J}$ be an
associated partition of unity,
and let $\varphi '_j\in C_0^\infty (U'_j)$ with $\varphi '_j=1$ on
$\operatorname{supp}\psi '_j$. By construction, $\psi '_1=1 $ on
$U_{10}$. We use the $\Phi _j$ and $\kappa _j$ on these subsets
(setting $\kappa _j(U'_j)=V'_j$, $\kappa _j(U_{10})=V_{10}$), and
denote the induced mappings for sections by $\Phi ^*_j$.

Now the auxiliary operator $P$ is taken to act as follows:
\begin{equation} 
Pu=\sum_{1\le j\le J}\varphi '_j[P_0((\psi '_ju)\circ {\Phi
^*_j}^{-1})]\circ \Phi _j^*.
\end{equation} 
 It is elliptic with positive definite principal symbol, and for
sections supported in $U_{10}$, it acts like $P_0$ when carried over
to $V_{10}$ (being a differential operator, it is local). The
resolvent $(P-\lambda )^{-1}$, defined for large $\lambda $ on the
rays in $\mathbb C\setminus \mathbb R_+$, is of course not local, 
but its symbol in the
local chart $V_1\times\mathbb C^M$ is, for $x\in V_{10}$,
equivalent with the symbol $(|\xi |^m+1-\lambda )^{-1}I_M$ of
$(P_0-\lambda )^{-1}$. For resolvents of 
differential operators, $q(x,\xi ,\lambda )\sim q_0(x,\xi ,\lambda )$
means that the difference is of order $-\infty $ and $O(\lambda
^{-N})$ for any $N$ (the symbols are strongly polyhomogeneous). 
This difference
does not affect the coefficient of $(-\lambda )^{-1}$ that we are
after.

Let $a(x,\xi )$ denote the symbol of $\psi _1A\varphi _1$ carried over to
$V_1\times 
\mathbb C^M$; it vanishes for $x\notin V_{10}$. Then the
symbol of $\psi _1A\varphi _1(P-\lambda )^{-1}$ on $V_1$ is equivalent with
$a(x,\xi )(|\xi |^m+1-\lambda )^{-1}$ (with an error that is
$O(\lambda ^{-N})$, any $N$); we use that the symbol composition here
gives only one (product) term.

To find the coefficient of $(-\lambda )^{-1}$ in the expansion of 
\linebreak $\Tr(\psi _1A\varphi _1(P-\lambda )^{-1})$, we now just have to
analyze the diagonal kernel calculated in $V_1$:
\begin{equation} 
\aligned K&(\psi _1A\varphi _1(P-\lambda )^{-1},x,x )\sim \int_{\mathbb R^n}  a(x,\xi )(|\xi
|^m+1-\lambda )^{-1}\,\dslash\xi \\
&\sim  
\sum_{ j\in \mathbb N  } \tilde c_{ j}(x)(-\lambda ) ^{\frac{\sigma +n -j}m-1}+ 
\sum_{k\in \mathbb N}\bigl( \tilde c'_{ k}(x)\log (-\lambda ) +\tilde c''_{
k}(x)\bigr)(-\lambda ) ^{{ -k}-1}.
\endaligned\label{3.2}\end{equation}
Here the value can be found explicitly, as follows.

Set $a_{-n}=0$ if $\sigma+n\notin \mathbb N$, and decompose the symbol
$a$ in three pieces $a_{>-n}$, 
$a_{-n}$ and $a_{<-n}$, where\begin{equation}
\aligned
a_{>-n}(x,\xi )&=\sum_{0\le j< \sigma +n}a_{\sigma -j}(x,\xi ),\\
a_{<-n}(x,\xi )&=a(x,\xi )-a_{-n}(x,\xi )-a_{>-n}(x,\xi ).
\endaligned
\label{3.3}\end{equation}
The symbol terms $a_{\sigma -j}(x,\xi )$ are assumed to be $C^\infty $ in
$(x,\xi )$  and homogeneous
of degree $\sigma -j$ in $\xi $ for $|\xi |\ge 1$. 
The strictly homogeneous version $a^h_{\sigma -j}$ is homogeneous for
$\xi \ne 0$ and coincides with $a_{\sigma -j}$ for $|\xi |\ge 1$.
For the terms in $a_{>-n}$, the strictly homogeneous versions are
integrable in $\xi $ at $\xi =0$. 

We recall that $\tslint a(x,\xi )\,\dslash\xi $ is defined for each $x$
as a finite part integral (in the sense of Hadamard),
namely the constant term in the asymptotic expansion 
of $\int_{|\xi |\le R}a(x,\xi )\,\dslash\xi $ in 
powers $R^{-m_j}$ and $\log R$, for $R\to\infty $. Here
\begin{equation}
\aligned
\slint a_{\sigma -j}(x,\xi )\,\dslash\xi &=
\int_{|\xi |\le 1}(a_{\sigma -j}(x,\xi
)-a^h_{\sigma -j}(x,\xi ))\,\dslash\xi, \text{ for }\sigma -j>-n,\\
\slint a_{-n}(x,\xi )\,\dslash\xi &=\int_{|\xi |\le 1}a_{-n}(x,\xi
)\,\dslash\xi,\\
\slint a_{<-n}(x,\xi )\,\dslash \xi&=\int_{\mathbb R^n}a_{<-n}(x,\xi
)\,\dslash \xi 
;\endaligned
\label{3.4}\end{equation}
as one can check using polar coordinates (the formulas are special
cases of [G4], (1.18)).

\begin{lemma}\label{Lemma 3.1} For $a_{\sigma -j}(x,\xi )$ with $\sigma
-j>-n$,\begin{equation} 
\aligned
\int_{\mathbb R^n}a_{\sigma -j}(x,\xi )(|\xi |^m-\lambda )^{-1}\,\dslash\xi 
&=(-\lambda
)^{\frac{\sigma +n-j}m-1}\int_{\mathbb R^n}a^h_{\sigma -j}(x,\xi )(|\xi
|^m+1)^{-1}\,\dslash\xi \\
&\quad+(-\lambda )^{-1}\slint a_{\sigma -j}(x,\xi
)\,\dslash\xi +O(\lambda ^{-2}),
\endaligned\label{3.5}\end{equation}
for $\lambda \to\infty $ on rays in $\mathbb C\setminus\mathbb R_+$. 
Also $\int_{\mathbb R^n}a_{\sigma -j}(x,\xi )(|\xi |^m+1-\lambda
)^{-1}\,\dslash\xi $ has an expansion in powers of $(-\lambda )$ plus
$o(\lambda ^{-1})$; here the coefficient of $(-\lambda )^{-1}$ is
likewise  $\tslint a_{\sigma -j}(x,\xi )\,\dslash\xi $.

For $a_{ -n}(x,\xi )$ one has:\begin{equation}
\aligned
\int_{\mathbb R^n}a_{-n}(x,\xi )(|\xi |^m-\lambda )^{-1}\,\dslash\xi 
&=\tfrac1m(-\lambda
)^{-1}\log(-\lambda )\int_{|\xi |=1}a_{-n}(x,\xi)\,\dslash S(\xi) \\
&\quad+(-\lambda )^{-1}\slint a_{-n}(x,\xi )\,\dslash\xi+O(\lambda ^{-2}),
\endaligned\label{3.10}\end{equation}
for $\lambda \to\infty $ on rays in $\mathbb C\setminus\mathbb R_+$. 
 $\int_{\mathbb R^n}a_{-n}(x,\xi )(|\xi |^m+1-\lambda
)^{-1}\,\dslash\xi $ has a similar expansion, the coefficient of
$(-\lambda )^{-1}$ again being $\tslint a_{-n}(x,\xi )\,\dslash\xi $.

\end{lemma}  

\begin{proof}
By homogeneity, we have for $\lambda \in\mathbb C\setminus{\overline{\mathbb R}_+}$, writing $\lambda =-|\lambda |e^{i\theta }$,
$|\theta |<\pi $,
\begin{equation}\aligned
\int_{\mathbb R^n}a^h_{\sigma -j}(x,\xi )&(|\xi |^m+|\lambda |e^{i\theta }
)^{-1}\,\dslash\xi \\
&=|\lambda |^{\frac{\sigma -j+n}m-1}\int_{\mathbb R^n}a^h_{\sigma -j}(x,\eta  )(|\eta
|^m+e^{i\theta }
)^{-1}
\,\dslash\eta .\endaligned\label{3.6 }\end{equation} 
This equals the first term in the right hand side of (\ref{3.5}) if $\theta
=0$, and the identity extends analytically to general $\lambda $.
Moreover, since \begin{equation}
\aligned
(|\xi |^m-\lambda )^{-1}&=
(-\lambda )^{-1}(1-|\xi |^m/\lambda )\\
&=(-\lambda )^{-1}\sum_{k\in\mathbb N}(|\xi
|^m/\lambda )^{k}\text{ for }|\lambda |\ge 2,\, |\xi |\le 1,
\endaligned\label{3.7
}\end{equation}
we find that for $|\lambda |\ge 2$,
\begin{equation}
\aligned
\int_{\mathbb R^n}&(a_{\sigma -j}(x,\xi )-a^h_{\sigma -j}(x,\xi ))(|\xi
|^m-\lambda )^{-1}\,\dslash\xi \\
&=\int_{|\xi |\le 1}(a_{\sigma -j}(x,\xi )-a^h_{\sigma -j}(x,\xi ))(|\xi
|^m-\lambda )^{-1}\,\dslash\xi \\
&=(-\lambda )^{-1}\int_{|\xi |\le 1}(a_{\sigma -j}(x,\xi )-a^h_{\sigma
-j}(x,\xi ))
\,\dslash\xi +O(\lambda ^{-2}).\endaligned\label{3.8}\end{equation}
This shows (\ref{3.5}), in view of (\ref{3.4}).

For the next observation, we use that the preceding results give an
expansion in powers of $1-\lambda $; then since\begin{equation}\aligned
(1-\lambda )^k &=\sum_{0\le l\le k}b_l\lambda ^l\text{ when
}k\in\mathbb N,\\
(1-\lambda )^s&=(-\lambda )^{s}+\sum_{ l\ge 1}b_l\lambda
^{s-l}\text{ when 
}|\lambda |\ge 2,\; s\notin \mathbb N,
\endaligned\label{3.9}\end{equation}
only $\tslint a_{\sigma -j}$ contributes to the coefficient of
$(-\lambda )^{-1}$.
 
Now consider $a_{-n}(x,\xi )$; again we can let $\theta =0$.
Here we write
\begin{equation}
\aligned
\int_{\mathbb R^n}&a_{-n}(x,\xi )(|\xi |^m-\lambda )^{-1}\,\dslash\xi 
\\&=\int_{|\xi |\ge 1 }a^h_{-n}(x,\xi )(|\xi |^m+|\lambda |)^{-1}\,\dslash\xi
+\int_{|\xi |\le 1}a_{-n}(x,\xi )(|\xi |^m-\lambda )^{-1}\,\dslash\xi.
\endaligned\end{equation}
The first term gives
\begin{equation}
\aligned
\int_{|\xi |\ge 1}&a^h_{-n}(x,\xi )(|\xi |^m+ |\lambda |)^{-1}\,\dslash\xi =
|\lambda | ^{-1}\int_{|\eta |\ge |\lambda |^{-1/m} }a^h_{-n}(x,\eta )(|\eta
|^m+1 )^{-1}\,\dslash\eta  \\
&= 
|\lambda |^{-1}\int_{r\ge |\lambda |^{-1/m} }r^{-1}(r^m+1 )^{-1}\,dr\int_{|\xi |=1}a_{-n}(x,\xi
)\,\dslash S(\xi)\\
&=\tfrac1m(-\lambda )^{-1}\log(-\lambda )\int_{|\xi |=1}a_{-n}(x,\xi
)\,\dslash S(\xi)+O(\lambda ^{-2}),
\endaligned\end{equation}
since $\int r^{-1}(r^m+1)^{-1}\,dr=\frac1m\int
s^{-1}(s+1)^{-1}\,ds=\frac1m(\log s-\log(s+1))$, $s=r^m$, where $\log(s+1)=O(s)$. The second term
gives, as in (\ref{3.8}),\begin{equation}
\int_{|\xi |\le 1}a_{-n}(x,\xi )(|\xi |^m-\lambda )^{-1}\,\dslash\xi =(-\lambda )^{-1}\int_{|\xi |\le 1}a_{-n}(x,\xi )
\,\dslash\xi +O(\lambda ^{-2}).
\end{equation} 
This implies (\ref{3.10}), in view of (\ref{3.4}).
The last statement follows using (\ref{3.9}).
 
\end{proof}

For $a_{<-n}$, it is very well known that\begin{equation}
\aligned
\int_{\mathbb R^n}a_{<-n}(x,\xi )(|\xi |^m-\lambda )^{-1}\,\dslash \xi &=
(-\lambda )^{-1}\int_{\mathbb R^n}a_{<-n}(x,\xi )\,\dslash \xi +o(\lambda
^{-1}),\\
&=
(-\lambda )^{-1}\slint a_{<-n}(x,\xi )\,\dslash \xi +o(\lambda
^{-1});
\endaligned\label{3.11}\end{equation}  
also here, 
\begin{equation}
\int_{\mathbb R^n}a_{<-n}(x,\xi )(|\xi |^m+1-\lambda )^{-1}\,\dslash \xi =
(-\lambda )^{-1}\slint a_{<-n}(x,\xi )\,\dslash \xi +o(\lambda
^{-1})\label{3.12}\end{equation}  
follows by use of (\ref{3.9}).
 
Collecting the informations, we have:

\begin{proposition}\label{Proposition 3.3} The coefficient of $(-\lambda )^{-1}$ in
the expansion {\rm (\ref{3.2})} for $K(\psi _1A\varphi _1(P-\lambda
)^{-1},x,x)$
equals
\begin{equation} 
\tilde c_n(x)+\tilde c''_0(x)=\slint  a(x,\xi )\,\dslash\xi .
\label{3.13
}\end{equation} 
\end{proposition} 

In the same localization, when we calculate $\psi _1A\varphi _1\log
P$ by a Cauchy 
integral (as in [G6]), the localized piece will give $\psi
_1A\varphi _1
\operatorname{OP}(\log (|\xi |^m+1))$. The symbol 
 of this operator is \begin{equation}
 r(x,\xi )=a(x,\xi )\log (|\xi |^m+1), 
\end{equation}
which has an expansion\begin{equation}
r(x,\xi )=a(x,\xi )(m\log |\xi |-|\xi
|^{-m}-\sum_{j\ge 2}d_j|\xi |^{-jm}) ,
\end{equation}
convergent for  $|\xi |\ge 2$.
Inserting the expansion of $a$ in homogeneous terms, we find since
$m>\sigma +n$ that the full
term of order $-n$ in $r(x,\xi )$ is $a_{-n}m\log |\xi |$,
with no log-free part. So 
\begin{equation}\operatorname{res}_{x,0}r=\int_{|\xi
|=1}\tr r_{-n,0}(x,\xi )\dslash S(\xi )=0.
\end{equation} 

It follows that the coefficient of
$(-\lambda )^{-1}$ in the trace expansion of $\psi _1A\varphi
_{1}(P-\lambda )^{-1}$ is
\begin{equation}
\int_{\mathbb R^n}\slint \tr a(x,\xi
)\,\dslash\xi dx 
=\int_{\mathbb R^n}(\TR_x
(\psi _1A\varphi _1)-\tfrac1m\operatorname{res}_{x,0}(\psi _1A\varphi
_1\log P))\,dx,
\label{3.14}\end{equation}
using that $\operatorname{res}_{x,0}(\psi _1A\varphi _1\log P)$ is 0. 
This shows formula (\ref{0.1}) in
this very particular case:
\begin{equation}
C_0(\psi _1A\varphi _1,P)=\int_{\mathbb R^n}(\TR_x 
(\psi _1A\varphi _1)-\tfrac1m\operatorname{res}_{x,0}(\psi _1A\varphi
_1\log P))\,dx.\label{3.15}
\end{equation}

Now, to find $C_0(\psi _1A\varphi _1,P')$ for a general operator $P'$
of order $m'\in\mathbb R_+$,
we combine (\ref{3.15}) with the trace 
defect formula (\ref{2.13}).
This gives, in the considered local coordinates: 
\begin{equation}
\aligned
C_0(\psi _1A\varphi _1,P')&=C_0(\psi _1A\varphi _1,P)+C_0(\psi
_1A\varphi _1,P')-C_0(\psi _1A\varphi _1,P)\\
&=\int_{\mathbb R^n}(\TR_x 
(\psi _1A\varphi _1)-\tfrac1m\operatorname{res}_{x,0}(\psi _1A\varphi
_1\log P))\,dx\\
&\quad- \operatorname{res}(\psi _1A\varphi _1(\tfrac1{m'}\log P'-\tfrac1m\log P))\\
&=\int_{\mathbb R^n}(\TR_x 
(\psi _1A\varphi _1)-\tfrac1{m'}\operatorname{res}_{x,0}(\psi
_1A\varphi _1\log P'))\,dx.
\endaligned
\end{equation}
To this we can add:
\begin{equation}
\aligned
C_0&(\psi _1A(1-\varphi _1),P')=\Tr(\psi _1A(1-\varphi _1))\\
&=\int(\TR_x 
(\psi _1A(1-\varphi _1))-\tfrac1{m'}\operatorname{res}_{x,0}(\psi
_1A(1-\varphi _1)\log P'))\,dx,
\endaligned
\end{equation}
where both terms have a meaning on $X$; $\TR_x$ defines the ordinary
trace integral and 
$\operatorname{res}_{x,0}$ is zero.

The method applies likewise to all the other terms $\psi _jA\varphi
_j$. Collecting the terms, and relabelling $P'$ of order $m'$ as $P$ of
order $m$, we have found:

\begin{theorem}\label{Theorem 3.4}
Let $A$ be a 
classical $\psi $do of order $\sigma \in\mathbb R$, and let $P$ be a classical elliptic $\psi $do's of
order $m\in\mathbb R_+$ such that the principal symbol has no eigenvalues on
$\mathbb R_-$. We have in local coordinates as used above:
\begin{equation}
\aligned
C_0(A,P)&=\sum_{1\le j\le J}C_0(\psi _jA,P), \text{ where}\\
C_0(\psi _jA,P)&=\int(\TR_x 
(\psi _jA)-\tfrac1{m}\operatorname{res}_{x,0}(\psi
_jA\log P))\,dx
\endaligned
\label{3.16}\end{equation}
(the contribution from $\psi _jA\varphi _j$ defined in the
corresponding local chart and that from $\psi _jA(1-\varphi _j)$
defined as an ordinary trace).
\end{theorem}  

Note that $C_0(A,P)$ is independent of how we localize, so the
expression resulting from (\ref{3.16}) is independent of the choice of
localization. 

The logarithm is here defined by cutting the complex plane along
$\mathbb R_-$. If $P$ is given with another ray free of eigenvalues, the
formulas hold with the logarithm defined to be cut along this ray.
(We do not here study the issue of how these expressions depend on
the ray.)

The invariance of the density
$\bigl(\TR_x(A)-\tfrac1m\operatorname{res}_{x,0}(A\log P)\bigr)\,dx$
in the formula (\ref{0.1}) is verified in [PS] also by a direct calculation.
We note that (\ref{0.1}) gives back the known formula
\begin{equation} 
C_0(A,P)=\TR(A)
\label{3.17}
\end{equation}
in cases where $\operatorname{res}_{x,0}(A\log P)$ vanishes. This is
so when $\sigma +n\notin\mathbb N$ ([KV], [L]), and also in cases
$\sigma +n\in\mathbb N$ with parity properties ([KV], [G4]): We say
that $A$ has even-even alternating parity (in short: is even-even), resp.\
has even-odd alternating parity (in short: is even-odd), when 
\begin{equation} \aligned
a_{\sigma -j}(x,-\xi)&=(-1)^{\sigma -j}a_{\sigma -j}(x,\xi ),\text{
resp.\ } \\
a_{\sigma -j}(x,-\xi)&=(-1)^{\sigma -j-1}a_{\sigma -j}(x,\xi
),
\endaligned
\label{3.18}
\end{equation}
for $|\xi |\ge 1$, all $j$. When $P$ is even-even of even order $m$, then the
classical part of $\log P$ is even-even. Then if (a) or (b) is satisfied:

(a) $A$ is even-even and $n$ is odd,

(b) $A$ is even-odd and $n$ is even,

\noindent $\operatorname{res}_{x,0}(A\log P)$ 
vanishes, $\TR_xA\,dx$ is a globally defined density, and
(\ref{3.17}) holds. [KV] treats the case (a), calling the even-even
operators odd-class (perhaps because they have a canonical trace in
odd dimension). The statements  
are extended to log-polyhomogeneous operators in [G4].
Observe the consequence:

\begin{corollary}\label{Corollary 4.1}
When {\rm (a)} or {\rm (b)} is satisfied,
$\operatorname{res}_{x,0}(A\log P)\,dx$ defines a global density for any 
$P$.
\end{corollary} 

\begin{proof} In these cases, since $\TR_x A\,dx$ defines a global density,
the other summand in
$\bigl(\TR_x(A)-\tfrac1m\operatorname{res}_{x,0}(A\log P)\bigr)\,dx$
must do so too.
\end{proof}

\appendix
\section{Corrections to earlier papers }
\medskip

{\it Correction to} [GS2]:
In Corollary 2.10 on page 45, the formulas in (2.38) for the
expansion coefficients $\tilde a_{j,l}$ are true only for $l=m_j$.
For $l<m_j$, the $\tilde a_{j,l}$ depend on the full set
$\{a_{j,l}\mid 0\le j\le m_j\}$. This is so, because the Taylor
expansion of $\Gamma (1-s)^{-1}$ must be taken into account.
\medskip

Hence in the comparison of the expansion of $\Tr(A(P-\lambda )^{-1})$
with $\zeta (A,P,s)$, only the primary coefficient at each pole of
$\Gamma (s)\zeta (A,P,s)$ is directly proportional to a coefficient
in $\Tr(A(P-\lambda )^{-1})$. Similar statements hold for comparisons
with $\Tr(A(P-\lambda )^{-N})$.
This has lead to systematic inaccuracies in a number
of subsequent works, however without substantial damage to the
results in general. 

We explain the needed correction in detail for [G4]
(where we have a few additional corrections), and then list the
related modifications needed in other papers.
\medskip

{\it Corrections to }[G4]:
\smallskip

1) The statements on page 69 linking the coefficients in (1.1) with
the coefficients in (1.2) with the same index by universal
proportionality factors is 
incorrect if $\nu +n\in\mathbb N$; the direct
proportionality holds only for the 
primary pole coefficients, not for the next Laurent coefficient at
each pole. Instead, at the second-order poles $-k, k\in\mathbb N$, there are
universal linear transition 
formulas linking the coefficient set for $(-\lambda )^{-k-N}\log(-\lambda
)$ and $(-\lambda )^{-k-N}$ with the coefficient set for $(s+k)^{-2}$
and $(s+k)^{-1}$. 

This follows from [GS2, Prop.\ 2.9], (3.21), as explained in
\ref{Lemma 1.1} of the present paper. The coefficients 
of $\Tr(A(P-\lambda )^{-N})$ at integer powers are directly
proportional to the Laurent 
coefficients of the
meromorphic function $\psi (s)$, where (with $N-1$ denoted $M$)\begin{equation}
\aligned\zeta (A,P,s)&=\tfrac{M!}{(s-M)\dots(s-1)}\tfrac1\pi \sin(\pi
(s-M))\psi (s-M),\\
\Gamma (s)\zeta (A,P,s)&=\tfrac{M!}{\Gamma (M-s)}\psi
(s-M).\endaligned\label{A.1}\end{equation} 
(Cf. (\ref{1.14}), use that $\frac1\pi \sin(\pi
(s-M))=(-1)^M/[\Gamma (s-1)\Gamma (1-s)]$.)
In calculations of Laurent series at the poles, the Taylor
expansion of the factor in front of $\psi (s)$ effects the higher terms.

Specifically in [G4], the sentence 
``The coefficients $\tilde c_j$ and $c_j$, $\tilde c'_k$ and $ c'_k$,
resp.\ $\tilde c''_k$ and $c''_k$ are proportional by universal nonzero
constants.'' should be replaced by: ``The coefficients $\tilde c_j$ and
$c_j$, resp.\ $\tilde c'_k$ and $ c'_k$,
are proportional by universal nonzero
constants. When the $c'_k$ vanish (e.g., when $\nu +n\notin \mathbb N$),
the same holds for $\tilde c''_k$ and $ c''_k$. More generally,
the pair $\{\tilde c'_k, \tilde c''_k\}$ is for each $k$
universally related to the pair $\{c'_k, c''_k\}$ in a linear way.'' The
statement ``$\tilde 
c''_0=c''_0$'' should be replaced by ``$\tilde 
c''_0=c''_0$ when $c'_0=0$'', and  
the description of $C_0(A,P)$ in terms of resolvent trace
expansion coefficients should be replaced by the description
given in the present paper in Section 1. 

However, since this changes the formula for $C_0(A,P)$ only by a
multiple of $\res 
A$, the results of [G4] on $C_0(A,P)$ remain valid, because they are
concerned with cases where $\res A=0$. The statements Th.\
1.3 (ii), Cor.\ 1.5 (ii) on the vanishing of all log-coefficients in
parity cases still imply the vanishing of all double poles in (1.2). 

In Section 3, the coefficients in (3.32) are linked with those in (3.30)
in a more complicated way than stated, where only the leading
coefficient at a pole is directly proportional to a coefficient in (3.30). 
But again, the results for parity cases remain valid since the
needed correction terms vanish in these cases.

\smallskip

2) Page 79, remove the factor 2 (twice) in formula (1.44).
\smallskip

3) Page 84, formulas (3.9) and (3.10): The sums over $k'$ should be
removed, and so should the additional term $-P^{-s-1}$ in the first
line. So $\mathcal P_l(P)=(-\log P)^l$ for all $l$.

\smallskip
4) Page 91, in formula (3.47), $-P^{-1}$ should be removed. 
\medskip

{\it Correction to }[G1]: 
Page 92, lines 7--8 from below, replace ``The coefficients in (9.10)
are proportional to those in 
(9.9) by universal factors.'' by ``The unprimed coefficients in
(9.10) are proportional to those in 
(9.9) by universal factors. For each $k$, the pair $\{\tilde a_{i,k},
\tilde a'_{i,k}\}$ (resp.\ $\{\tilde b_{i,k},
\tilde b'_{i,k}\}$) is universally related to the pair $\{ a_{i,k},
 a'_{i,k}\}$ (resp.\ $\{b_{i,k},
b'_{i,k}\}$) in a linear way.''

\medskip

{\it Corrections to }[G2]:
Page  4, lines 5--7 from below, replace ``The coefficients $\tilde
a_k$, $\tilde a'_k$ and $\tilde a''_k$ are proportional to the
coefficients $a_k$, $a'_k$ and $ a''_k$ in (0.1) (respectively)
by universal nonzero proportionality factor (depending on $r$).''
 should be replaced by: ``The coefficients $\tilde
a_k$ and $\tilde a'_k$ are proportional to the
coefficients $a_k$ and $a'_k$ in (0.1) (respectively)
by universal nonzero proportionality factor (depending on $r$).
For each $k\ge 0$, the pair $\{\tilde a'_k, \tilde a''_k\}$ is
universally related to the pair $\{a'_k, a''_k\}$ in a linear way.'' 
\medskip

{\it Corrections to }[G3]:
Page 262, lines 13--14 should be changed as for [G2] above. In
line 15, ``$a''_0(F)=\tilde a''_0(F)$'' should be replaced by:
``$a''_0(F)=\tilde a''_0(F)$ if $\tilde a'_0(F)=0$''. There are some
consequential reformulations in Sections 4 and 5, which do not
endanger the results since the $a''_0$ terms are characterized in
general modulo local contributions (and $a'_0$ is such), with precise
statements only when $a'_0=0$.
\medskip

{\it Corrections to }[G5]: The statement on Page 44, lines
11--12 from below ``There are some 
universal proportionality factors linking the coefficients $\tilde
c_j$ and $c_j$, $\tilde c'_k$ and $ c'_k$, 
resp.\ $\tilde c''_k$ and $c''_k$'' should be changed as indicated
for [G4].

\medskip

{\it Correction to }[GSc1]:
Page 171, line 1, ``The coefficients $\tilde c_j$,
$\tilde c'_l$ , $\tilde c''_l$ are proportional to the coefficients 
$ c_j$, $ c'_l$ , $ c''_l$ by universal
constants.'' should be replaced by ``The coefficients $\tilde c_j$,
$\tilde c'_l$  are proportional to the coefficients 
$ c_j$, $ c'_l$ by universal
constants.''  

\medskip

{\it Corrections to }[GSc2]:
The definition of $C_0(A,P)$ on page 1644 and the statements in
lines 6--8 on page 1645 should be modified as for [G4]. This has
no consequences for the results, which are mainly concerned with the trace
definitions {\it modulo local contributions}, with exact formulas established
only when the residue corrections vanish.

\medskip

{\it Corrections to }[G6]:
The order $m$ of $P_1$ and $P_2$ should be taken to be even in Theorem
3.10 and in the second half of Theorem 4.5. This is in order to
assure that the classical parts of the $\log P_i$ satisfy the transmission
condition. For $m$ odd, the local formulas for the coefficients in terms of
integrals with $\log\lambda $ are still valid, but the interpretation as
residues is not covered by [FGLS]. Similar remarks for Theorem 5.2.


\begin{thebibliography}{0}

\bibitem[FGLS]{} 
B. V. Fedosov, F. Golse, E. Leichtnam and E. Schrohe,
{\it The noncommutative residue for manifolds with boundary}
 J. Funct. Anal.
{\bf 142}, 
1--31 (1996).


\bibitem[G1]{} 
G. Grubb, {\it Trace expansions for 
pseudodifferential boundary problems for Dirac-type operators
and more general systems},  
Arkiv f. Mat.
{\bf 37}, 45--86, (1999).




\bibitem[G2]{}
G. Grubb,
{\it Logarithmic terms in trace expansions of Atiyah-Patodi-Singer
problems}, 
  Ann. Global Anal. Geom.
{\bf  24},   1--51
(2003).

 

\bibitem[G3]{} G. Grubb,
 {\it  Spectral boundary conditions for generalizations of Laplace and Dirac operators},
  Comm. Math. Phys.
{\bf  240}, 
  243--280 (2003).

 

\bibitem[G4]{} G. Grubb,
{\it  A resolvent approach to traces and zeta Laurent expansions},
Contemp.\ Math.\ 
{\bf  366},   67--93
(2005). Corrected in  arXiv: math.AP/0311081. 
 

\bibitem[G5]{} G. Grubb,
{\it  Analysis of invariants associated with spectral boundary
problems for elliptic operators}, 
Contemp.\
Math.\ {\bf  366},
  43--64 (2005).
 
\bibitem[G6]{}
G. Grubb, 
{\it  On the logarithm component in trace defect formulas}, 
Comm.\ Part.\ Diff.\ Equ.\ {\bf 30}, 1671--1716 (2005).

 

\bibitem[GH]{} G. Grubb and L. Hansen
{\it  Complex powers of resolvents of pseudodifferential operators},
  Comm. Part. Diff. Equ.
{\bf  27}, 
  2333--2361
(2002).
 


\bibitem[GSc1]{}  G. Grubb and E. Schrohe,
{\it  Trace expansions and the noncommutative residue for manifolds
with boundary},
  J. Reine Angew. Math.
{\bf  536},
  167--207 (2001).
 

\bibitem[GSc2]{} G. Grubb and E. Schrohe,
{\it  Traces and quasi-traces on the Boutet de Monvel algebra}
  Ann. Inst. Fourier
{\bf  54},
  1641--1696
(2004).
 

\bibitem[GS1]{} G. Grubb and R. Seeley, 
{\it  Weakly parametric    pseudodifferential
    operators and Atiyah-Patodi-Singer boundary problems},
      Invent. Math. {\bf  121},  481--529 (1995). 


\bibitem[GS2]{}G. Grubb and R. Seeley, 
{\it  Zeta and eta functions for Atiyah-Patodi-Singer
     operators},  J. Geom.\ Anal.\ {\bf  6},   31--77 (1996).
 



\bibitem[Gu]{} V. Guillemin,
{\it  A new proof of Weyl's formula on the asymptotic distribution
of eigenvalues},  Adv. Math. {\bf  102},
  184--201 (1985).
 

\bibitem[KV]{} M. Kontsevich and S. Vishik,
{\it  Geometry of determinants of elliptic operators},
Functional Analysis on the Eve of the 21'st Century, Vol. I (New
Brunswick, N.J.\ 1993),
 Progr. Math. 131, Birkh\"auser, Boston,
  173--197 (1995).

 
\bibitem[L]{}
 M. Lesch, {\it  On the noncommutative residue for
pseudodifferential operators with log-poly\-ho\-mo\-ge\-ne\-ous symbols},
  Ann. Global Anal. Geom.
{\bf  17},  151--187 (1999).

 
\bibitem[Lo]{}
 P. Loya,
{\it The structure of the resolvent of elliptic pseudodifferential
operators}, 
J. Funct. Anal.
{\bf  184},   77--134 (2001).
 

\bibitem[MN]{}
 R. Melrose and V. Nistor,
{\it  Homology of pseudodifferential operators I. Manifolds with
boundary}, arXiv: funct-an/9606005.
 

\bibitem[O]{} K. Okikiolu,
{\it  The multiplicative anomaly for determinants of elliptic
operators},  Duke Math. J.
{\bf  79},  723--750 (1995).
 

\bibitem[PS]{}
 S. Paycha and S. Scott,
{\it  An explicit Laurent expansion for regularized integrals of
holomorphic symbols}, to appear, arXiv: math.AP/0506211.
 

\bibitem
[W]{} M. Wodzicki,
{\it  Local invariants of spectral asymmetry},  Invent. Math.
{\bf 75},  143--178 (1984).

 
\end{thebibliography}
\end{document}